\documentclass[12pt,reqno]{amsart}

\usepackage{amscd,amssymb,amsthm,amstext,hyperref,txfonts}
\usepackage{rotating,graphicx,graphics,subfigure,epsfig,psfrag}
\usepackage{geometry} 

\usepackage[parfill]{parskip}    

\usepackage{empheq,verbatim} 

\usepackage[all,knot]{xy}
\xyoption{arc}

\def\proof{{\it Proof) }} 
\def\qed{$\Box$}

\def\FF{{\mathcal F}}
\def\HH{{\mathcal H}}

\def\NN{{\mathcal N}}

\DeclareSymbolFont{AMSb}{U}{msb}{m}{n}
\DeclareMathSymbol{\N}{\mathbin}{AMSb}{"4E}
\DeclareMathSymbol{\Z}{\mathbin}{AMSb}{"5A}
\DeclareMathSymbol{\R}{\mathbin}{AMSb}{"52}
\DeclareMathSymbol{\Q}{\mathbin}{AMSb}{"51}
\DeclareMathSymbol{\I}{\mathbin}{AMSb}{"49}
\DeclareMathSymbol{\C}{\mathbin}{AMSb}{"43}

\theoremstyle{plain}
\newtheorem{theorem}{Theorem}[section]

\newtheorem{lemma}[theorem]{Lemma}

\newtheorem{corollary}[theorem]{Corollary}

\newtheorem{question}[theorem]{Question}

\theoremstyle{definition}
\newtheorem{definition}[theorem]{Definition}

\newtheorem{remark}[theorem]{Remark}

\renewcommand{\labelenumi}{(\arabic{enumi})}

\def\lll{\langle}
\def\rrr{\rangle}

\def\link{\textrm{link}}
\def\pmo{{}^{-1}}

\newtheorem*{theoremnone}{Theorem}

\newcommand{\con}{\mathrm{CO}}
\newcommand{\cob}{\overline{\mathrm{CO}}}
\newcommand{\eqag}{=_{A(\Gamma)}}

\date{\today}
\author[Sang-hyun Kim]{Sang-hyun Kim}
\address{}
\email{}
\title{Co-contractions of Graphs and Right-angled Artin Groups}

\begin{document}
\begin{abstract}
We define an operation on finite graphs, called {\em co-contraction}.
Then we show that for any co-contraction $\hat\Gamma$ of a finite graph
$\Gamma$,
the right-angled Artin group on $\Gamma$
 contains 
a subgroup which is isomorphic to the right-angled Artin group
on $\hat\Gamma$.
As a corollary, 
we exhibit a family of graphs, without any induced cycle of length at least 5,
such that the right-angled Artin groups on those graphs contain hyperbolic
surface groups. This gives the negative answer to a question raised by
Gordon, Long and Reid.
\end{abstract}
\maketitle
%
%
%
%

\section{Introduction}\label{sec:intro}
In this paper, by a {\em graph} we mean a finite graph without loops
and without multi-edges. 
A {\em right-angled Artin group} is a group 
defined by a presentation with a finite
generating set,
where the relators are certain commutators between the generators.
Such a presentation naturally determines the {\em underlying graph},
where the vertices correspond to the generators and the edges to
the pairs of commuting generators. It is known that the isomorphism type
of a right-angled Artin group uniquely
 determines the isomorphism
type of the underlying graph~\cite{droms1987,KMNR1980}. 
Also, right-angled Artin groups
possess various group theoretic properties.
To name a few, right-angled
Artin groups 
are 
linear~\cite{humphries1994,HW1999,DJ2000},
biorderable~\cite{DT1992},
biautomatic~\cite{vanwyk1994}
and  moreover,
admitting free and cocompact actions
 on finite-dimensional CAT(0) cube complexes~\cite{CD1995,MV1995,NR1998}.

On the other hand,
it is interesting to ask what 
we can say about the isomorphism type
of the underlying graph,
if a right-angled Artin group satisfies a given 
group theoretic property.
Let $\Gamma$ be a graph. We denote the vertex set and 
the edge set of $\Gamma$ by $V(\Gamma)$ and $E(\Gamma)$, respectively.
The {\em complement graph}  of $\Gamma$ is the graph $\overline\Gamma$
defined 
by
$
V(\overline\Gamma)=V(\Gamma)$ and 
$E(\overline\Gamma)=\{ \{u,v\}: \{u,v\}\not\in
E(\Gamma)\}$.
For a subset $S$ of $V(\Gamma)$ the
 {\em induced subgraph}  
 on $S$, denoted by $\Gamma_S$, is defined to be
the maximal subgraph of $\Gamma$ with the vertex set $S$. 
This implies that
$V(\Gamma_S)=S$ and 
$E(\Gamma_S)=\{\{u,v\}:u,v\in S\textrm{ and }\{u,v\}\in E(\Gamma)\}$.
If $\Lambda$ is another graph, 
an {\em induced} $\Lambda$ in $\Gamma$  means 
an induced subgraph isomorphic to $\Lambda$ in $\Gamma$. 
$C_n$ denotes the cycle of length
$n$. That is, $V(C_n)$ is a set of $n$ vertices, say $\{v_1,v_2,\ldots,
v_n\}$, and $E(C_n)$ consists of the edges $\{v_i,v_j\}$ where
$|i-j|\equiv 1\ (\textrm{mod } n)$.
Let $A(\Gamma)$ be the right-angled Artin group with its
underlying graph $\Gamma$. Then, the following are true. 
\begin{itemize}
\item
$A(\Gamma)$ is 
coherent, if and only if 
$\Gamma$ is {\em chordal}, i.e. $\Gamma$ does not contain
an induced $C_n$ for any $n\ge4$~\cite{droms1987a}.
This happens if and only if 
$[A(\Gamma),A(\Gamma)]$ is free~\cite{SDS1989}.
\item
$A(\Gamma)$ is 
a virtually 3-manifold group, if and only if each connected component 
of $\Gamma$ is tree or triangle~\cite{droms1987a,gordon2004}
\item
$A(\Gamma)$ is 
subgroup separable, if and only if no induced subgraph of $\Gamma$
is a square or a path of length 3~\cite{MR2006}. This happens
if and only if every subgroup of $A(\Gamma)$ is also a right-angled
Artin group~\cite{droms1987b}.
\item
$A(\Gamma)$
 contains a {\em hyperbolic
surface group}, i.e.
the fundamental group
of a closed, hyperbolic surface,
 if there exists an
induced $C_n$ for some $n\ge5$ in $\Gamma$~\cite{SDS1989,CW2004}. 
\end{itemize} 

In~\cite{GLR2004}, Gordon, Long and Reid proved that 
a word-hyperbolic (not necessarily right-angled) Coxeter group
either is virtually-free or contains a hyperbolic surface group.
They also showed that certain (again, not necessarily right-angled)
Artin groups do not contain
a hyperbolic surface group, raising the following question.

\begin{question}  \label{que:main}
Does $A(\Gamma)$ contain a hyperbolic surface group
if and only if $\Gamma$ contains an induced $C_n$ for some $n\ge5$?
\end{question}

In this paper, we give the negative answer to the above question.
 Let
$\Gamma$ be a graph and $B$ be a set of vertices of $\Gamma$ such 
that $\Gamma_B$ is connected. The
{\em  contraction} of $\Gamma$ relative to $B$ is 
the graph $\con(\Gamma,B)$  obtained from $\Gamma$
by collapsing $\Gamma_B$ to a vertex, and deleting
loops or multi-edges. We define
the {\em co-contraction}  $\cob(\Gamma,B)$ 
of  
$\Gamma$ relative $B$, such that
$\cob(\Gamma,B)=\overline{\con(\overline{\Gamma},B)}$. 
Then we prove
 the following theorem, which will imply that
 $A(\overline{C_n})$ contains $A(\overline{C_5})=A(C_5)$ and hence
 a hyperbolic surface subgroup, for $n\ge5$ (see 
Figure~\ref{fig:barc6contract}). An easy combinatorial
 argument shows that  $\overline{C_n}$ 
 does not
contain an induced cycle of length at least 5, for $n>5$. 
\begin{theoremnone} 
Let $\Gamma$ be a graph and $B$ be a set of vertices in $\Gamma$,
such that $\overline{\Gamma_B}$ is
connected. 
Then $A(\Gamma)$ contains a subgroup isomorphic to
$A(\cob(\Gamma,B))$.
\end{theoremnone}

In this paper, the above theorem
is proved in the following steps.

In Section~\ref{sec:prelim}, we recall basic 
facts on right-angled Artin groups and 
HNN extensions. A {\em dual van Kampen diagram} is
described.
We owe the notations to~\cite{CW2004} where a closely
related concept, a {\em dissection}, was defined and
used with great clarity.

In Section~\ref{sec:contraction}, we define  {\em co-contraction}
of a graph, and examine its properties.

In Section~\ref{sec:contractanticnt}, we prove the theorem by 
exhibiting an embedding
of $A(\cob(\Gamma,B))$ into $A(\Gamma)$. The main tool for the proof
is a dual van Kampen diagram.

In Section~\ref{sec:contraction_words}, we compute intersections of
certain subgroups of right-angled Artin groups.
From this, we  deduce a more detailed version of the theorem
describing some other choices of the embeddings.

{\bf Acknowledgement.} I am 
deeply grateful to my thesis advisor, Andrew Casson, 
for
his insights and guidance. I would like to also thank Daniel Spielman, for 
helpful comments.

%
%
%
%

\section{Preliminary on Right-angled Artin Groups}  \label{sec:prelim} 
Let $\Gamma$ be a graph. 
The {\em right-angled
Artin group on $\Gamma$} is the group presented as,
\[ A(\Gamma)=\lll v\in V(\Gamma)\;|\;[a,b]=1\textrm{ if and only if }
\{a,b\}\in E(\Gamma)\rrr\]
Each element of $A(\Gamma)$ can be expressed as
$w=\prod_{i=1}^k  c_i^{e_i}$, where $ c_i\in V(\Gamma)$ and
$e_i=\pm1$. 
Such an expression is called a {\em word (of length $k$)} and 
each $ c_i^{e_i}$ is called a {\em letter} of the word $w$.
We say the word $w$ is 
{\em reduced}, if the length is minimal among the words
representing the same element. 
For each $i_0=1,2,\ldots,k$, the word $w_1=\prod_{i={i_0}}^k c_i^{e_i} \cdot
\prod_{i=1}^{i_0-1} c_i^{e_i}$ is called a {\em cyclic conjugation} of 
$w=\prod_{i=1}^k  c_i^{e_i}$.
By a {\em subword} of $w$, we mean a word 
$w'=\prod_{i=i_0}^{i_1} c_i^{e_i}$ for some $1\le i_0<i_1\le k$.
A letter or a subword $w'$ of $w$ is on the {\em left}
of a letter or a subword $w''$ of $w$, if $w'=\prod_{i=i_0}^{i_1} c_i^{e_i}$
and $w''=\prod_{i=j_0}^{j_1} c_i^{e_i}$ for some $i_1<j_0$.

The expression $w_1=w_2$ 
shall mean that $w_1$ and $w_2$ are equal as words
(letter by letter). On the other hand,
$w_1=_{A(\Gamma)}w_2$ means that the words
$w_1$ and $w_2$ represent the same element in $A(\Gamma)$. 
For an element
$g\in A(\Gamma)$ and a word $w$, 
$w\eqag g$ means that the word $w$ is representing the group element
$g$. $1$ denotes both the trivial element in $A(\Gamma)$ and
the empty word, depending on the context.

\begin{figure}[htb!] 
\psfrag{w0}{{$w=c^{-1}aba^{-1}b^{-1}c$}}
\psfrag{T}{{$\downarrow$}}
\psfrag{a}{{$a$}}
\psfrag{b}{{$b$}}
\psfrag{c}{{$c$}}
\psfrag{am}{{$a^{-1}$}}
\psfrag{bm}{{$b^{-1}$}}
\psfrag{cm}{{$c^{-1}$}}
\psfrag{(a)}{{(a) $\tilde\Delta$}}
\psfrag{(b)}{{(b) $(\tilde\Delta)^*$}}
\psfrag{(c)}{{(c) $(\tilde\Delta)^*\setminus B(v_\infty)$}}
\psfrag{vinf}{{$v_\infty$}}
\psfrag{Bv}{{$B(v_\infty)$}}
\centerline{\epsfig{file=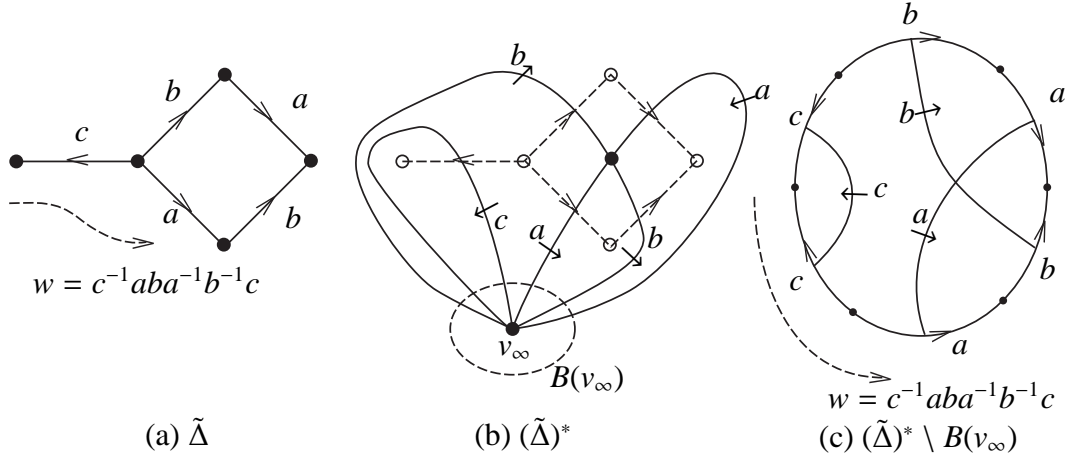, width=.95\hsize}}
\caption{Constructing a dual van Kampen diagram from 
a van Kampen diagram, for $w=c^{-1}aba^{-1}b^{-1}$ in 
$\lll a,b,c\;|\;[a,b]=1\rrr$. 
\label{fig:dualvk}}
\end{figure}

Let $w$ be a word representing the trivial element 
in $A(\Gamma)$. 
A {\em     dual van Kampen diagram $\Delta$} for $w$ in $A(\Gamma)$ 
is a pair $(\HH,\lambda)$
satisfying the following (Figure~\ref{fig:dualvk} (c)):
\renewcommand{\labelenumi}{(\roman{enumi})}
\begin{enumerate}
\item
$\HH$ is a set of transversely oriented simple closed curves and 
transversely oriented 
properly embedded arcs in general position, 
in an oriented disk $D\subseteq \R^2$.
\item
$\lambda$ is a map from $\HH$
to 
$V(\Gamma)$ such that
$\gamma$ and $\gamma'$ in $\HH$ are intersecting
{\em only if} $\lambda(\gamma)$ and $\lambda(\gamma')$
 are adjacent in $\Gamma$.
\item
Enumerate the boundary points of the arcs in $\HH$ as $v_1,v_2,\ldots,
v_m$ so that $v_i$ and $v_j$ are adjacent on $\partial D$ 
if and only if $|i-j|\equiv 1\ (\textrm{mod } n)$.
For each $i$, let $a_i$ be the label of the arc that
intersects with $v_i$. Put $e_i=1$ if, at $v_i$, the orientation of 
$\partial D$ coincides with the transverse orientation of the
arc that $v_i$ is intersecting, and $e_i=-1$ otherwise. 
Then $w$ is a cyclic conjugation of $v_1^{e_1} v_2^{e_2}  
\cdots v_m^{e_m}$.
\end{enumerate}

Note that simple closed curves in 
a dual van Kampen diagram can always be assumed to be removed. Also,
we may assume that two curves in $\Delta$ are minimally intersecting,
in the sense that there does not exist any bigon formed by arcs in $\HH$.
See~\cite{CW2004} for more details, as well as generalization of this definition
to arbitrary compact surfaces, rather than a disk.

Let $\tilde\Delta\subseteq S^2$ be  
a (standard) van Kampen diagram 
for $w$, with respect to a standard presentation 
$A(\Gamma)=\lll V(\Gamma)\;|\; [u,v]=1 \textrm{ if and only if }
\{u,v
\}\in E(\Gamma)\rrr$
(Figure~\ref{fig:dualvk}).
Consider $\tilde\Delta^*$, the dual of $\tilde\Delta$ in $S^2$,
and name the vertex which is dual to the face $S^2\setminus \tilde
\Delta$ as $v_\infty$.
Then for a sufficiently small ball $B(v_\infty)$ around $v_\infty$,
$\tilde\Delta^*\setminus B(v_\infty)$
can be considered as a dual van Kampen diagram
with a suitable choice of the labeling map.
Therefore a dual van Kampen diagram exists for
any word $w$ representing the trivial element in 
$A(\Gamma)$. 
Conversely,
a van Kampen diagram $\tilde\Delta$ for a word
can be obtained from 
a dual van Kampen diagram  $\Delta$
by considering
the dual complex again. So, the existence of a dual van Kampen
diagram for a word $w$ implies that $w\eqag 1$.

Given a dual van Kampen diagram $\Delta$, divide $\partial D$ 
into segments so that each segment intersects with exactly one
arc in $\HH$. Let the label and the orientation of each segment be
induced from  those of the arc that intersects with the segment. The
resulting labeled
and directed graph on $\partial D$ is called the {\em boundary} of $\Delta$
and denoted by $\partial \Delta$.

We call each arc in $\HH$ labeled by $q\in V(\Gamma)$  as 
a {\em $q$-arc}, and each segment in $\partial \Delta$ labeled by 
$q$ as a {\em $q$-segment}. Sometimes we identify
the letter $q^{\pm1}$ 
of $w$
 with the corresponding $q$-segment.
 A connected union of segments on $\partial\Delta$
is called an {\em interval}.
By convention, a subword $w_1$ of $w$ shall also denote
the corresponding interval (called {\em $w_1$-interval})
on $\partial\Delta$.

Now let  
$\Delta=(\HH,\lambda)$
be
 a dual van Kampen diagram
  on $D\subseteq\R^2$.
  Suppose
$\gamma$ is a  properly embedded arc in $D$,
 which is either an element in $\HH$ or 
in general position with
$\HH$.
Then one can cut 
$\Delta$ along $\gamma$ in the
following sense.
First, cut $D$ along $\gamma$ to get two disks  $D'$ and $D''$.
Consider the intersections of the disks with 
the curves in $\HH$.
Then,
let those curves in $D'$ and $D''$ 
inherit the transverse orientations and the labeling maps from $\Delta$.
We obtain 
two dual van Kampen diagrams, one for each of $D'$ and $D''$.
Conversely, we can glue two dual van Kampen diagrams along identical
words.
An {\em innermost $q$-arc} $\gamma$
is
a $q$-arc such that the interior of $D'$ or $D''$ does not intersect
any $q$-arc.

\renewcommand{\labelenumi}{(\roman{enumi})}
\begin{definition} \label{def:qpair}
Let $\Gamma$ be a graph.
Let $w$ be a word representing the trivial element in $A(\Gamma)$,
and $\Delta$ be a dual van Kampen diagram for $w$. 
Two segments on the boundary of $\Delta$ are called 
a {\em cancelling $q$-pair}
if there exists a $q$-arc joining the segments. 
For {\em any} word
$w_1$,
two letters of $w_1$ are called a 
{\em cancelling $q$-pair}
if there exist another word $w_1'=_{A(\Gamma)}w_1$ 
and 
a dual van Kampen diagram $\Delta$ for $w_1w_1'^{-1}$, such that
the two letters are a $q$-pair with respect to $\Delta$.
A cancelling $q$-pair is also called as a {\em $q$-pair} for abbreviation.
A {\em cancelling pair} is a cancelling $q$-pair for some $q\in V(\Gamma)$.
\end{definition}

\renewcommand{\labelenumi}{(\arabic{enumi})}

For a group $G$ and its subset $P$, $\lll P\rrr$ denotes
the subgroup generated by $P$.  For a subgroup $H$ of $A(\Gamma)$,
$w\in H$ shall mean that $w$ represents an element in $H$.

\begin{lemma} \label{lem:qpair}
Let $\Gamma$ be a graph and $q$ be a vertex of $\Gamma$. 
If a word $w$ in $A(\Gamma)$ has a $q$-pair,
then
$w=w_1 q^{\pm1} w_2 q^{\mp1} w_3$ 
for some subwords $w_1,w_2$ and $w_3$ such that
$w_2\in\lll \link_{\Gamma}(q)\rrr$.
In this case,
$w$ is not reduced.
\end{lemma}

\proof

There exists a word $w'=_{A(\Gamma)}w$ 
and a dual van Kampen 
diagram $\Delta$ for $ww'^{-1}$, such that
a $q$-arc joins two segments of $w$.

\begin{figure}[htb!]
\psfrag{w0}{\footnotesize $w$}
\psfrag{w1}{\footnotesize $w_1$}
\psfrag{w2}{\footnotesize $w_2$}
\psfrag{w3}{\footnotesize $w_3$}
\psfrag{wp}{\footnotesize $w'$}
\psfrag{wt}{\footnotesize $\tilde w_2$}
\psfrag{g1}{\footnotesize $\gamma$}
\psfrag{g2}{\footnotesize $\Delta_0$}
\psfrag{g3}{\footnotesize $\beta'$}
\psfrag{bj}{\footnotesize $d_j$}
\psfrag{dj}{\footnotesize $\delta_j$}
\psfrag{del}{\footnotesize $\Delta$}
\psfrag{bpj}{\footnotesize $d_j^{f_j}$}
\psfrag{bmj}{\footnotesize $d_j^{-f_j}$}
\psfrag{bm}{\footnotesize $b^{-1}$}
\psfrag{b}{\footnotesize $b$}
\psfrag{amp}{\footnotesize $a^{\mp1}$}
\psfrag{apm}{\footnotesize $a^{\pm1}$}
\psfrag{q}{\footnotesize $q$}
\psfrag{qm}{\footnotesize $q^{\mp1}$}
\psfrag{bj}{\footnotesize $d_j$}
\psfrag{T}{$\downarrow$}
\centerline{\epsfig{file=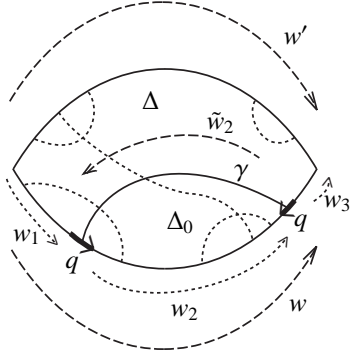, width=.30\hsize}}
\caption{Cutting $\Delta$ along $\gamma$. \label{fig:qpair}}
\end{figure}

Write $w=w_1 q^{\pm1} w_2 q^{\mp1} w_3$, where
the letters $q^{\pm1}$ and $q^{\mp1}$
(identified with the corresponding segments on 
$\partial \Delta$)  are joined by 
a $q$-arc $\gamma$ as in Figure~\ref{fig:qpair}.

Cut $\Delta$ along $\gamma$, to get a dual van Kampen diagram $\Delta_0$,
which contains $w_2$ on its boundary. 
Give $\Delta_0$ the orientation
that  coincides with  the orientation of $\Delta$ on $w_2$.
Let $\tilde w_2$ be the word,
read off by following $\gamma$ in the orientation of
 $\Delta_0$.
$\tilde w_2\in\lll\link_\Gamma(q)\rrr$,
for the arcs intersecting with $\gamma$ are labeld by 
vertices in 
$\link_\Gamma(q)$.
Since $\Delta_0$ is a dual van Kampen diagram for the word 
$w_2 \tilde w_2$, we have
$w_2=_{A(\Gamma)}\tilde w_2^{-1}\in
\lll\link_\Gamma(q)\rrr$.
\qed

For $S\subseteq V(\Gamma)$, we let $S^{-1}=\{q^{-1}:q\in S\}$ and 
$S^{\pm1}=S\cup S\pmo$.
The following lemma is standard, and we 
briefly sketch the proof. 

\begin{lemma} \label{lem:subgraph}
Let $\Gamma$ be a graph and $S$ be a subset of $V(\Gamma)$.
Then the following are true.
\begin{enumerate}
\item
$\lll S\rrr$ is 
isomorphic
to $A(\Gamma_S)$.
\item
Each letter of any reduced word in $\lll S\rrr$ 
is in $S^{\pm1}$.
\end{enumerate}
\end{lemma}

\proof 

(1)
The inclusion $V(\Gamma_S)\subseteq V(\Gamma)$ induces
a map $f:A(\Gamma_S)\rightarrow A(\Gamma)$. Let $w$ be a word
representing an element in $\ker f$. Since $w=_{A(\Gamma)}1$,
there exists a 
dual van Kampen diagram $\Delta$ for the word $w$ in $A(\Gamma)$.
Remove simple closed curves
labeled by $V(\Gamma)\setminus V(\Gamma_S)$, if there is any.
Since the boundary of $\Delta$ is labeled by vertices in
 $V(S)$,
$\Delta$ can be considered as a dual van Kampen diagram for 
a word $w$ in $A(\Gamma_S)$.
So we get 
$w=_{A(\Gamma_S)}1$. 

(2)

$w=_{A(\Gamma)}w'$ for some word $w'$ such that the letters of $w'$
are in $S$. 
Let $\Delta$ be a dual van Kampen
diagram for $w w'\pmo$. If $w$ contains a $q$-segment
for some $q\not\in S$, then a $q$-arc joins two segments
in $\Delta$, and these segments must be in $w$.
This is impossible by
Lemma~\ref{lem:qpair}. 
\qed

From this point on, $A(\Gamma_S)$ is considered
as a subgroup of $A(\Gamma)$, 
for $S\subseteq V(\Gamma)$.
Let $H$ be a group and $\phi:C\rightarrow D$ be an isomorphism
between subgroups of $H$. 
Then we define
$H\ast_\phi
=\lll H,t\;|\;t\pmo c t = \phi(c),\textrm{ for }c\in C\rrr$,
which is the HNN extension of $H$ 
with the amalgamating map $\phi$ and the stable letter $t$.
Sometimes, we explicitly state what the stable letter is.
If $C=D$ and $\phi$ is the identity map, then we let
$H\ast_C
=\lll H,t\;|\; t\pmo c t = t\textrm{ for }c\in C\rrr$.

For a vertex $v$ of a graph $\Gamma$,
the {\em link of $v$}
is the set
\[\link_\Gamma(v)=\{ u\in V(\Gamma) : u\textrm{ is adjacent to }v\}\]

\begin{lemma} \label{lem:artinhnn}
Let $\Gamma$ be a graph.
Suppose $\Gamma'$ is an  induced subgraph of $\Gamma$ such that
$V(\Gamma')=V(\Gamma)\setminus\{v\}$ for some $v\in V(\Gamma)$.
Let 
$C$ be the subgroup of $A(\Gamma')$ generated by
$\textrm{link}_\Gamma
(v)
$.
Then the inclusion $A(\Gamma')\hookrightarrow A(\Gamma)$
extends to the isomorphism
$f:A(\Gamma')\ast_C
\rightarrow
A(\Gamma)$ such that $f(t)=v$.
\end{lemma}

\proof Immediate from the definition of right-angled Artin groups.\qed

We first note the following general lemma.

\begin{lemma} \label{lem:jkt}
Let $H$ be a group and
 $\phi:C\rightarrow D$ be an isomorphism between
 subgroups $C$ and $D$.
Suppose $K$ is a subgroup of $H$ and
$J=\lll K,t\rrr\le H\ast_\phi$. We let
$\psi:J\cap C\rightarrow J\cap D$ be the restriction 
of $\phi$.
Then
the inclusion $J\cap H\hookrightarrow J$
extends to the  isomorphism 
$f:(J\cap H)
\ast_\psi
\rightarrow
J$
 such that $f(\hat t)=t$, where $\hat t$ and $ t$ denote 
 the stable letters of $(J\cap H)\ast_\psi$ and
 $H\ast_\phi$,
 respectively.
\end{lemma}

\proof
Note that $G=H\ast_\phi$ acts on a tree $T$, 
with a vertex $v_0$ and an
 edge $e_0=\{v_0,t.v_0\}$
satisfying 
$\textrm{Stab}( v_0) = H$ and $\textrm{Stab} (e_0) = C$ 
~\cite{serre2003}.
Let $T_0$ be the induced subgraph on $\{j.v_0:j\in J\}$.
For each vertex $j.v_0$ of $T_0$, write $j=
k_1 t^{\epsilon_1}
k_2 t^{\epsilon_2}
\cdots
k_m t^{\epsilon_m}$, where $k_i\in K$ and $\epsilon_i=\pm1$ for each $i$.
Then the following sequence in $V(T_0)$
\begin{eqnarray*}
v_0 & = &k_1. v_0,\\
k_1 t^{\epsilon_1}. v_0 & = &
k_1 t^{\epsilon_1} k_2. v_0,\\
k_1 t^{\epsilon_1} k_2 t^{\epsilon_2}.  v_0 & = &
k_1 t^{\epsilon_1} k_2 t^{\epsilon_2} k_3.  v_0,\\
\ldots\\
k_1 t^{\epsilon_1} k_2 t^{\epsilon_2} k_3 \cdots t^{\epsilon_m}.  v_0 &=& j.v_0
\end{eqnarray*}
gives rise to a path in $T_0$ from $v_0$ to $j.v_0$. Hence
$T_0$ is connected. Note that 
$\psi:J\cap C=\textrm{Stab}_J(e_0)\rightarrow
J\cap D=\textrm{Stab}_J (e_0)^t$.
Since 
$J$ acts on a tree $T_0$, 
we have an isomorphism
$J\cong \textrm{Stab}_J(v_0)\ast_\psi
=
(J\cap H)\ast_\psi$.\qed

%
\section{Co-contraction of Graphs}\label{sec:contraction}
%
%
Let $\Gamma$ be a graph and $B\subseteq V(\Gamma)$.
We say  $B$ is  {\em connected}, if $\Gamma_B$ is connected.
$B$ is  {\em anticonnected}, if $\overline{\Gamma_B}$ is connected.

\renewcommand{\labelenumi}{(\roman{enumi})}

\begin{definition}
Let $\Gamma$ be a graph and $B\subseteq V(\Gamma)$.
\begin{enumerate}
\item
If $B$ is connected, the
{\em contraction of $\Gamma$ relative to $B$} is the graph
$\con(\Gamma,B)$ defined by:
\begin{eqnarray*}
V(\con(\Gamma,B))& = &\left(V(\Gamma)\setminus B\right)\cup\{v_B\}\\
E(\con(\Gamma,B))& = &E(\Gamma_{V(\Gamma)\setminus B})
\cup\{\{v_B,q\}:q\in V(\Gamma)\setminus B\textrm{ and }
\link_\Gamma(q)\cap B\ne\varnothing\}
\end{eqnarray*}
\item
If $B$ is anticonnected, the 
 {\em co-contraction of $\Gamma$ relative to $B$} is the graph
$\cob(\Gamma,B)$ defined by:
\begin{eqnarray*}
V(\cob(\Gamma,B))& = &\left(V(\Gamma)\setminus B\right)\cup\{v_B\}\\
E(\cob(\Gamma,B))& = &E(\Gamma_{V(\Gamma)\setminus B})
\cup\{\{v_B,q\}:q\in V(\Gamma)\setminus B\textrm{ and }\link_\Gamma(q)\supseteq B\}
\end{eqnarray*}
\item
More generally, if $B_1,B_2,\ldots,B_m$ are disjoint connected 
subsets of $V(\Gamma)$, then  inductively define
\[
\con(\Gamma,(B_1,B_2,\ldots,B_m))= \con(
\con(\Gamma,(B_1,B_2,\ldots,B_{m-1})),B_m)\]
and if $B_1,B_2,\ldots,B_m$ are disjoint anticonnected
 subsets, then similarly,
\[
\cob(\Gamma,(B_1,B_2,\ldots,B_m))= \cob(\cob(\Gamma,(B_1,B_2,\ldots,B_{m-1})),B_m)\]
\end{enumerate}
\end{definition}

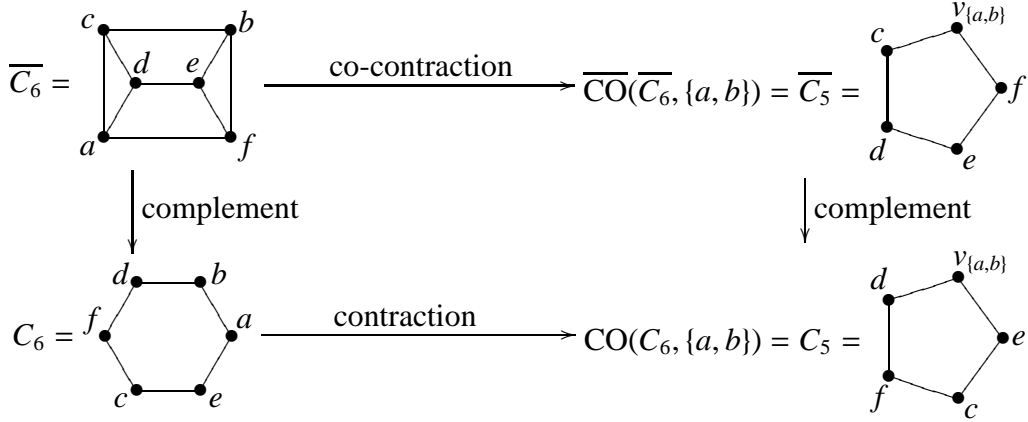
\begin{figure}[ht!] 
\[
\xymatrix{
\overline{C_6}=
\xy
0;/r.10pc/:
(20, 17)*{}="p";
(25, 20)*{b};
(-20, 17)*{}="a";
(-25, 20)*{c};
(-10, 0)*{}="b";
(-8, 7)*{d};
(-20,-17)*{}="c";
(-25, -20)*{a};
(20, -17)*{}="r";
(25, -20)*{f};
(10, 0)*{}="q";
(8, 6)*{e};
"p"*{\bullet};
"a"*{\bullet};
"b"*{\bullet};
"c"*{\bullet};
"r"*{\bullet};
"q"*{\bullet};
"p";"q"**\dir{-};
"p";"r"**\dir{-};
"q";"r"**\dir{-};
"a";"b"**\dir{-};
"a";"c"**\dir{-};
"b";"c"**\dir{-};
"a";"p"**\dir{-};
"b";"q"**\dir{-};
"c";"r"**\dir{-};
\endxy
\ar[rrrr]^*{\mbox{co-contraction}\qquad\quad}
\ar[d]^{\mbox{complement}}
& & & &
\cob(\overline{C_6},\{a,b\}) = 
\overline{C_5}
=
\xy
0;/r.10pc/:
(22, 20)*{}="a"; 
(29, 25)*{v_{\{a,b\}}};
(0, 13)*{}="b"; 
(-3,18)*{c};
(0, -11)*{}="c"; 
(-3,-18)*{d};
(22, -18)*{}="d"; 
(26,-22)*{e};
(36, 1)*{}="e";
(41, 1)*{f};
"a"*{\bullet};
"b"*{\bullet};
"c"*{\bullet};
"d"*{\bullet};
"e"*{\bullet};
"a";"b"**\dir{-};
"b";"c"**\dir{-};
"c";"d"**\dir{-};
"d";"e"**\dir{-};
"e";"a"**\dir{-};
\endxy
\ar[d]^{\mbox{complement}}
\\
C_6 =
\xy
0;/r.10pc/:
(10, 17)*{}="p";
(16, 20)*{b};
(-10, 17)*{}="a";
(-15, 20)*{d};
(-20, 0)*{}="b";
(-24, 5)*{f};
(-10,-17)*{}="c";
(-15, -20)*{c};
(10, -17)*{}="r";
(15, -20)*{e};
(20, 0)*{}="q";
(24, 5)*{a};
"p"*{\bullet};
"a"*{\bullet};
"b"*{\bullet};
"c"*{\bullet};
"r"*{\bullet};
"q"*{\bullet};
"p";"q"**\dir{-};
"q";"r"**\dir{-};
"a";"b"**\dir{-};
"b";"c"**\dir{-};
"a";"p"**\dir{-};
"c";"r"**\dir{-};
\endxy
\ar[rrrr]^*{\mbox{contraction}\qquad\qquad}
&&&& \con(C_6,\{a,b\})=C_5
=
\xy
0;/r.10pc/:
(22, 20)*{}="a"; 
(29, 25)*{v_{\{a,b\}}};
(0, 13)*{}="b"; 
(-3,20)*{d};
(0, -11)*{}="c"; 
(-3,-18)*{f};
(22, -18)*{}="d"; 
(26,-22)*{c};
(36, 1)*{}="e";
(41, 1)*{e};
"a"*{\bullet};
"b"*{\bullet};
"c"*{\bullet};
"d"*{\bullet};
"e"*{\bullet};
"a";"b"**\dir{-};
"b";"c"**\dir{-};
"c";"d"**\dir{-};
"d";"e"**\dir{-};
"e";"a"**\dir{-};
\endxy
}
\]
\caption{
Note that $\{q:\link_\Gamma(q)\supseteq \{a,b\} \}= 
\{c,f\}$, i.e. $c$ and $f$ are common neighbors of 
$\{a,b\}$. Hence
in $\cob(\overline{C_6},\{a,b\})$, 
$v_{\{a,b\}}$ is adjacent to
$c$ and $f$. This can be
also viewed by
looking at the complement graph of $\overline{C_6}$, namely $C_6$, and collapsing
the edge $\{a,b\}$.
\label{fig:barc6contract}}
\end{figure}

In a graph $\Gamma$,
if $B$ is connected, then $\con(\Gamma,B)$ is obtained by
(homotopically) 
collapsing $\Gamma_B$ 
onto one vertex and removing any loops or multi-edges. 
If $B$ is anticonnected,
one has (see Figure~\ref{fig:barc6contract})
\[ \cob(\Gamma,B)=\overline{\con(\overline\Gamma,B)}\]
 If $B\subseteq V(\Gamma)$ and
$\link_\Gamma(q)\supseteq B$,
then we say that $q$ is a {\em common neighbor of} $B$.

The following lemma
states that the co-contraction of  a set of
anticonnected
vertices can be
obtained by considering a sequence of co-contractions of two 
non-adjacent
vertices. The proof is immediate by considering the complement graphs.

\begin{lemma} \label{lem:contractinduct}
Let $\Gamma$ be a graph and $B\subseteq V(\Gamma)$ be anticonnected.
Then there exists a sequence of graphs
\[\Gamma_0=\Gamma,\Gamma_1,\Gamma_2,\ldots,\Gamma_p=\cob(\Gamma,B)\]
such that for each $i=0,1,\ldots,p-1$,
$\Gamma_{i+1}$ is a co-contraction of $\Gamma_i$ relative
to a pair of non-adjacent vertices of $\Gamma_i$.\qed
\end{lemma}

\begin{lemma} \label{lem:barcn}
\begin{enumerate}
\item
If $B$ is a connected subset of $p$ vertices of $C_n$, then
$\con(C_n,B)\cong C_{n-p+1}$.
\item
If $B$ is an anticonnected subset of $p$ vertices of $\overline{C_n}$, then
$\cob(\overline{C_n},B)\cong \overline{C_{n-p+1}}$.
\end{enumerate}
\end{lemma}

\proof (1) is obvious. Considering
the complement graphs, (2) follows from (1).\qed

%
\section{Co-contraction of Graphs and Right-angled 
Artin Groups}\label{sec:contractanticnt}
%
%

Let $\Gamma$ be a graph and $B$ be an anticonnected subset of $V(\Gamma)$.
Fix a word $\tilde w\in\lll B\rrr$ in $A(\Gamma)$. 
If a vertex $x$ of $\cob(\Gamma,B)$ is adjacent to $v_B$,
then
$x$ is a common neighbor of $B$ in $\Gamma$, and so,
$[\phi(x),\phi(v_B)]=[x,\tilde w]\eqag1$. This implies that there exists a map 
$\phi:
A(\cob(\Gamma,B))
\rightarrow A(\Gamma)$ satisfying
\[\phi(x) =
\left\{
\begin{array}{ll}
\tilde w\quad & \textrm{if }x = v_{B}\\
x\quad& \textrm{if }x\in 
V(\cob(\Gamma,B))\setminus\{v_B\}=V(\Gamma)\setminus B 
\end{array}
\right.
\]
In this section,
we show that this map $\phi$ is injective for 
a suitable choice of the word $\tilde w$.
First, we prove the injectivity for
the case when $B=\{a,b\}$ and $\tilde  w=b\pmo a b$.

\renewcommand{\labelenumi}{(\roman{enumi})}

\begin{lemma} \label{lem:contract2}
Let $\Gamma$ be a graph. Suppose $a$ and $b$ are
non-adjacent vertices of $\Gamma$.
Then there exists an injective map
$
\phi:A(\cob(\Gamma,\{a,b\}))\rightarrow A(\Gamma)$
satisfying 
\[\phi(x) =
\left\{
\begin{array}{ll}
b\pmo a b\quad & \textrm{if }x = v_{\{a,b\}}\\
x\quad& \textrm{if } x\in V(\Gamma)\setminus\{a,b\}\\
\end{array}
\right.
\]
\end{lemma}

\proof

Let 
$\hat \Gamma=\cob(\Gamma,\{a,b\}), \hat v=v_{\{a,b\}}$ and
$A=\{q:q\in V(\Gamma)\setminus\{a,b\}\}$. For $q\in A$,
let $\hat q$ denote the corresponding vertex in $\hat\Gamma$, and 
$\hat A=\{\hat q:q\in A\}$.

Define $\phi:A(\hat\Gamma)\rightarrow A(\Gamma)$ by
\[\phi(x) =
\left\{
\begin{array}{ll}
b\pmo a b\quad &\textrm{if } x = \hat v\\
q\quad &\textrm{if } x=\hat q\in\hat A
\end{array}
\right.
\]

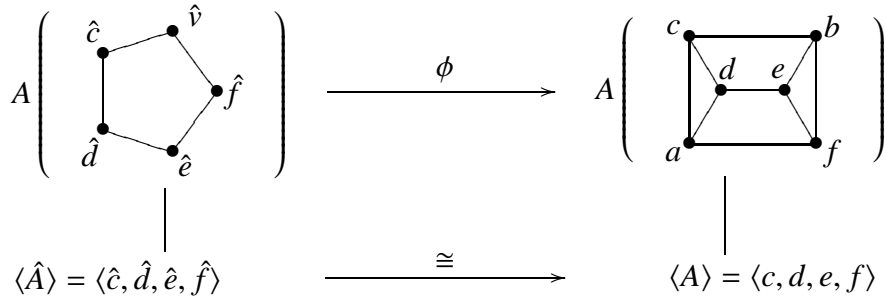
\begin{figure}[ht!]
\[
\xymatrix
{
A\left(\quad
\xy
0;/r.10pc/:
(22, 20)*{}="a"; 
(29, 25)*{\hat v};
(0, 13)*{}="b"; 
(-3,20)*{\hat c};
(0, -11)*{}="c"; 
(-4,-18)*{\hat d};
(22, -18)*{}="d"; 
(26,-22)*{\hat e};
(36, 1)*{}="e";
(41, 1)*{\hat f};
"a"*{\bullet};
"b"*{\bullet};
"c"*{\bullet};
"d"*{\bullet};
"e"*{\bullet};
"a";"b"**\dir{-};
"b";"c"**\dir{-};
"c";"d"**\dir{-};
"d";"e"**\dir{-};
"e";"a"**\dir{-};
\endxy\quad
\right)\quad
\ar[rrr]^*{\phi}
\ar@{-}[d]
& & &  
\quad A\left(\quad
\xy
0;/r.10pc/:
(20, 17)*{}="p";
(25, 20)*{b};
(-20, 17)*{}="a";
(-25, 20)*{c};
(-10, 0)*{}="b";
(-8, 7)*{d};
(-20,-17)*{}="c";
(-25, -20)*{a};
(20, -17)*{}="r";
(25, -20)*{f};
(10, 0)*{}="q";
(8, 6)*{e};
"p"*{\bullet};
"a"*{\bullet};
"b"*{\bullet};
"c"*{\bullet};
"r"*{\bullet};
"q"*{\bullet};
"p";"q"**\dir{-};
"p";"r"**\dir{-};
"q";"r"**\dir{-};
"a";"b"**\dir{-};
"a";"c"**\dir{-};
"b";"c"**\dir{-};
"a";"p"**\dir{-};
"b";"q"**\dir{-};
"c";"r"**\dir{-};
\endxy
\quad\right)
\ar@{-}[d]
\\
\lll \hat A\rrr = \lll \hat c,\hat d,\hat e,\hat f\rrr\qquad\quad
\ar[rrr]^*{\cong}
&&& 
\qquad\quad\lll A\rrr
=\lll c,d,e,f\rrr
}
\]
\caption{
An example of a co-contraction induced map between right-angled Artin groups.
\label{fig:cocontinj}}
\end{figure}

Suppose $\phi$ is not injective. Choose a word $\hat w$ of the minimal length
in $\ker \phi\setminus\{1\}$. 
Write $\hat w=\prod_{i=1}^k\hat c_i^{e_i}$, where $\hat c_i\in \hat A\cup
\{\hat v\}$
and
$e_i=\pm1$. Since $\hat\Gamma_{\hat A}$ is isomorphic to
$\Gamma_A$, 
$\phi$ maps $\lll\hat A\rrr$ isomorphically onto $\lll A\rrr$ (Figure~\ref{fig:cocontinj}).
So $\hat c_i = \hat v$ for some $i$.

\begin{figure}[htb!]
\psfrag{w}{\footnotesize 
$w=b^{-1}a^{\pm1}b\cdot w_1\cdot b^{-1}a^{\mp1}b\cdot w_2$}
\psfrag{w1}{\footnotesize $w_1$}
\psfrag{w2}{\footnotesize $w_2$}
\psfrag{a}{\footnotesize $a$}
\psfrag{al}{\footnotesize $\alpha$}
\psfrag{be}{\footnotesize $\beta$}
\psfrag{w_1'}{\footnotesize $w_1'$}
\psfrag{bm}{\footnotesize $b^{-1}$}
\psfrag{b}{\footnotesize $b$}
\psfrag{amp}{\footnotesize $a^{\mp1}$}
\psfrag{apm}{\footnotesize $a^{\pm1}$}
\psfrag{d0}{\footnotesize $\Delta_0$}
\centerline{\epsfig{file=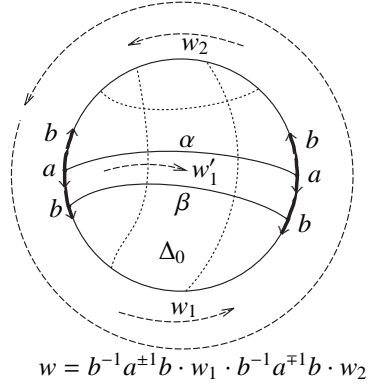, width=.30\hsize}}
\caption{$\Delta$ in the proof of Lemma~\ref{lem:contract2}.\label{fig:contract2}}
\end{figure}

Let $w=\prod_{i=1}^k \phi(\hat c_i)^{e_i}$. Since $w\eqag1$, 
there exists a  
dual van Kampen
diagram $\Delta=(\HH,\lambda)$ for $w$ in $A(\Gamma)$.
In $\Delta$, choose an innermost $a$-arc $\alpha$. 
By considering a cyclic conjugation of $\hat w$ if necessary, 
one may write
$\hat w = \hat v^{\pm1}\cdot \hat w_1 \cdot \hat v^{\mp1}\cdot \hat w_2$
and
$w = b^{-1}a^{\pm1}b\cdot w_1 \cdot b^{-1} a^{\mp1}b\cdot w_2$,
so that  $w_1=\phi(\hat w_1),w_2=\phi(\hat w_2)$ and 
$\alpha$ joins the  leftmost $a^{\pm1}$ of $w$
and
the $a^{\mp1}$  between $w_1$ and $w_2$ (Figure~\ref{fig:contract2}). Then the interval
$w_1$ does not contain any
$a$-segment.
Since each $b$-segment in $w$ is adjacent to 
some $a$-segment, one sees that
there does not exist any $b$-segment in $w_1$, either.
Hence,
$w_1\in\lll A\rrr= A(\Gamma_A)$ and $\hat w_1\in\lll\hat A\rrr
= A(\hat\Gamma_{\hat A})$. 
Note that $\Gamma_A\cong\hat\Gamma_{\hat A}$.
Since $\hat w_1$
is reduced, so is $w_1$.

Let $\beta$ be the $b$-arc that meets the letter 
$b$, following $a^{\pm1}$ on the left of $w_1$ in $w$.
$\beta$ does not intersect $\alpha$, for  $[a,b]\ne1$. 
Since $w_1$ does not contain any 
$b$-segment, $\beta$ intersects 
with
the letter $b^{-1}$ between $w_1$ and $w_2$.

$w_1$ does not contain any cancelling
pair, for $w_1$ is reduced. 
So each segment of $w_1$ is joined to a segment in $w_2$ by an arc in
$\HH$. Such an arc must intersect both $\alpha$ and $\beta$.
This implies that the segments in $w_1$ are labeled
by vertices in
$\link_\Gamma(a)\cap\link_\Gamma(b)
=\phi(\link_{\hat\Gamma}(\hat v))$.
It follows that $\hat w_1\in\lll
\link_{\hat\Gamma}(\hat v)\rrr$, from the following diagram.
\[
\xymatrix{
\ \ \qquad\hat w_1\in\lll\hat A\rrr\ar@{|->}[d] & 
&\lll\link_{\hat\Gamma}(\hat v)\rrr
\ar[d]^{\cong} &\le&
\lll\hat A\rrr\ar[d]_{\phi}^{\cong}\\
w_1 & \in & \lll\link_\Gamma(a)\cap\link_\Gamma(b)\rrr & \le &\lll A\rrr
}
\]
But then,
$\hat w=\hat v^{\pm1} \hat w_1\hat v^{\mp1}\hat w_2
=_{A(\hat\Gamma)}\hat w_1 \hat w_2$, which contradicts to the minimality of the 
length of $\hat w$.
\qed

\renewcommand{\labelenumi}{(\arabic{enumi})}

\begin{theorem} \label{thm:anticnt}
Let $\Gamma$ be a graph and $B$ be an anticonnected subset
of $V(\Gamma)$.
Then $A(\Gamma)$ contains a subgroup isomorphic to
$A(\cob(\Gamma,B))$
\end{theorem}

\proof
Proof is immediate from 
Lemma~\ref{lem:contractinduct} and Lemma~\ref{lem:contract2}.
\qed

Figure~\ref{fig:cocontinj}
and
Lemma~\ref{lem:contract2} show the existence of an isomorphism
\[
\phi: A(C_5)\rightarrow\lll b\pmo a b,c,d,e,f\rrr\le A(\overline{C_6})\]
More generally, we have the following corollary.

\begin{corollary} \label{cor:artincnbar}
\begin{enumerate}
\item
$A(\overline{C_n})$ contains a subgroup isomorphic to
$A(\overline{C_{n-p+1}})$ for each $1\le p\le n$.
\item
If $\Gamma$ contains an induced $C_n$ or $\overline{C_n}$ for some
$n\ge 5$, then $A(\Gamma)$ contains a hyperbolic surface group.
\end{enumerate}
\end{corollary}

\proof
(1)
Immediate from Lemma~\ref{lem:barcn} and Theorem~\ref{thm:anticnt}.

(2)  $A(C_n)$ contains a hyperbolic surface group for
$n\ge5$~\cite{SDS1989}. One has an embedding
$\phi:A(C_5)=A(\overline{C_5})\hookrightarrow A(\overline{C_n})$, for
$n\ge5$.\qed

A simple combinatorial argument shows that for $n> 5$, 
the induced subgraph of $\overline{C_n}$ 
on any five vertices contains a triangle. So 
$\overline{C_n}$ does not
contain an induced $C_m$ for any $m\ge 5$. From
 the Corollary~\ref{cor:artincnbar} (2),
we deduce the negative answer to  Question~\ref{que:main} as follows.

\renewcommand{\labelenumi}{(\roman{enumi})}
\begin{corollary} \label{cor:main}
There exists an infinite family $\FF$ of 
graphs satisfying the following.
\begin{enumerate}
\item
each element in $\FF$ does not contain an induced $C_n$ for $n\ge5$,
\item
each element in $\FF$ is  not 
an induced subgraph of another element in $\FF$,
\item
for each $\Gamma\in\FF$, $A(\Gamma)$ contains
a hyperbolic surface group.
\end{enumerate}
\end{corollary}

\proof Let $\FF=\{\overline{C_n}:n>5\}$.\qed
\renewcommand{\labelenumi}{(\arabic{enumi})}

A graph $\Gamma$ is called {\em weakly chordal} if 
$\Gamma$ does not contain an induced $C_n$ or $\overline{C_n}$ for
any
$n\ge 5$~\cite{hayward1985}.  Let $\NN=\{\Gamma: A(\Gamma)\textrm{ does not contain a hyperbolic
surface group}\}$. Corollary~\ref{cor:artincnbar}
shows that every graph in $\NN$
is weakly chordal.
Also,
Theorem~\ref{thm:anticnt} implies that $\NN$ is closed under
co-contraction.
On the other hand, if a graph is weakly chordal, 
then a co-contraction of the graph is also weakly chordal 
~\cite{kim2006}. This raises the following question.

\begin{question}
Does $A(\Gamma)$ contain a hyperbolic surface group if and only if
$\Gamma$ is weakly chordal?
\end{question}

\section{Contraction Words} \label{sec:contraction_words}

In Lemma~\ref{lem:contract2}, the word $b\pmo a b$ was used to construct
an injective map from $A(\cob(\Gamma,\{a,b\}))$ into $A(\Gamma)$.
This can be generalized by considering a {\em contraction word}, 
defined as 
follows.

\begin{definition} \label{def:contraction}
\begin{enumerate}
\item
Let $\Gamma_0$ be an anticonnected graph. A sequence
$b_1,b_2,\ldots,b_p$ of vertices of $\Gamma_0$
 is a {\em contraction sequence of
$\Gamma_0$}, if the following holds:
for any $(b,b')\in V(\Gamma_0)\times V(\Gamma_0)$,
there exists $l\ge1$ and $1\le k_1< k_2 <\cdots<k_l\le p$ such that,
$b_{k_1},b_{k_2},\ldots,b_{k_l}$ is a path from $b$ to $b'$ in $\overline\Gamma$.
\item
Let $\Gamma$ be a graph and $B$ be an anticonnected set of vertices of $\Gamma$.
A reduced word $w=\prod_{i=1}^p b_p^{e_i}$
is called a {\em contraction word of $B$} if 
$b_i\in B$, $e_i=\pm1$ 
for each $i$,  and
$b_1,b_2,\ldots,b_p$ is a contraction sequence of $\Gamma_B$.
An element of $A(\Gamma)$ is called a {\em contraction element},
if it can be represented by a contraction word.
\end{enumerate}
\end{definition}

\begin{remark}
If $a$ and $b$ are non-adjacent vertices in $\Gamma$,
then
any word in $\lll a,b\rrr\setminus\{a^mb^n:m,n\in\Z\}^{\pm1}$
is a contraction word of $\{a,b\}$.
\end{remark}

We first note the following general lemma.

\begin{lemma} \label{lem:conjugatepower}
Let $\Gamma$ be a graph and $g\in A(\Gamma)$.
Then $g=_{A(\Gamma)}u^{-1}vu$ for some words $u,v$ such that
$u\pmo v^m u$ is reduced for each $m\ne0$.
\end{lemma}

\proof
Choose words $u,v$ such that $u\pmo vu$ is a
reduced word representing $g$ and the length of $u$ is maximal. 
We will show that $u\pmo v^m u$ is reduced for any $m\ne 0$.

\begin{figure}[ht!]
\psfrag{w}{\footnotesize $w$}
\psfrag{u1}{\footnotesize $u_1$}
\psfrag{u2}{\footnotesize $u_2$}
\psfrag{u1m}{\footnotesize $u_1^{-1}$}
\psfrag{u2m}{\footnotesize $u_2^{-1}$}
\psfrag{v1}{\footnotesize $v_1$}
\psfrag{v2}{\footnotesize $v_2$}
\psfrag{v3}{\footnotesize $v_3$}
\psfrag{vi}{\footnotesize $v_i$}
\psfrag{ga}{\footnotesize $\gamma$}
\psfrag{vi1}{\footnotesize $v_{i+1}$}
\psfrag{vm}{\footnotesize $v_m$}
\psfrag{dd}{\footnotesize $\ldots$}
\psfrag{q}{\footnotesize $q$}
\psfrag{u}{\footnotesize $u$}
\psfrag{um}{\footnotesize $u^{-1}$}
\psfrag{w1}{\footnotesize $w_1$}
\psfrag{w0}{\footnotesize $w_0$}
\psfrag{w2}{\footnotesize $w_2$}
\psfrag{w3}{\footnotesize $w_3$}
\psfrag{(a)}{(a)}
\psfrag{(b)}{(b) Case 1.}
\psfrag{(c)}{(c) Case 2.}
\centerline{\epsfig{file=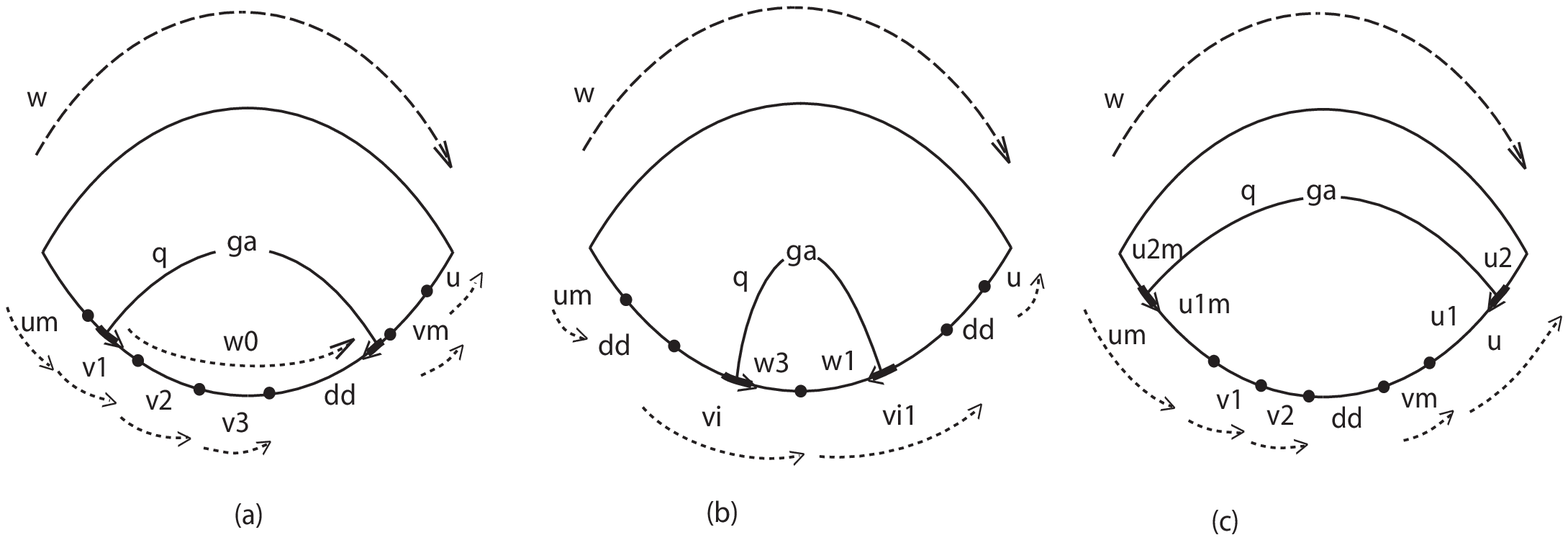, width=.95\hsize}}
\caption{Proof of Lemma~\ref{lem:conjugatepower}.
\label{fig:conjugatepower}}
\end{figure}

Assume that $u\pmo v^m u$ is 
not reduced for some $m\ne 0$. We may assume that $m>0$.
 Let $w$ be a reduced word for $u\pmo v^m u$.
Draw  a 
dual van Kampen diagram $\Delta$ for $u\pmo v^m u w^{-1}$.
Let $v_i$ denote the $v$-interval on
$\partial\Delta$ corresponding to the $i$-th occurrence of $v$ from 
the left in $u\pmo v^m u$ (Figure~\ref{fig:conjugatepower} (a)).

By Lemma~\ref{lem:qpair}, there exists a $q$-arc $\gamma$ joining
two $q$-segments of $u\pmo v^m u$ for  some $q\in V(\Gamma)$.
Let $w_0$ denote the interval between those two $q$-segments. 
We may choose $q$ and $\gamma$ so that the number of the segments
in $w_0$ is minimal. Then 
any arc intersecting with a segment in $w_0$ 
must intersect $\gamma$. It follows that any letter in $w_0$
should commute with $q$. Moreover, 
$w_0$ does not contain any 
$q$-segment.

{\em Case 1. The intervals $u\pmo$ and $u$ do not intersect with
$\gamma$.}

 Since $w_0$ does not contain any $q$-segment, 
$\gamma$ joins $v_i$ and $v_{i+1}$ for some $i$  (Figure~\ref{fig:conjugatepower} (b)).
Then one can write
$v=w_1 q^{\pm1} w_2 q^{\mp1} w_3$ for some subwords $w_1,w_2,w_3$ of $v$
such that
and $w_0=w_3 w_1$.
$[w_3,q]\eqag 1\eqag[w_1,q]$. 
So $u\pmo v u \eqag 
u\pmo q^{\pm1} w_1 w_2 w_3 q^{\mp1} u$,
which contradicts to the maximality of $u$.

{\em Case 2. $\gamma$ intersects $u$- or $u\pmo$-interval.}

Suppose $u\pmo$ intersects $\gamma$. Since $u\pmo v$ is reduced,
$\gamma$ cannot intersect $v_1$. So, $w_0$ contains $v_1$.
Since $w_0$ does not contain
any $q$-segment, $v$ does not contain the letters
$q$ or $q\pmo$ and so,
 $\gamma$ cannot intersect any 
$v_i$ for  $i=1,\ldots,m$.  
$\gamma$ should intersect with
the $u$-interval of $u\pmo v^m u$  (Figure~\ref{fig:conjugatepower} (c)).
This implies that
 $\gamma$ intersects with the leftmost $q$-segment in
the $u$-interval of $u\pmo v^m u$.
One can write
$u\pmo v^m u = u_2\pmo q^{\pm1} u_1\pmo v^m u_1 q^{\mp1} u_2$
such that any letter in $w_0=u_1\pmo v^m u_1$ commutes
with $q$, i.e.
$[q,u_1]\eqag1\eqag[q,v]$. But then
$u\pmo v u =_{A(\Gamma)} u_2\pmo u_1\pmo v u_1 u_2$,
which is a contradiction to the assumption that 
$u\pmo v u$ is reduced.\qed

\begin{lemma} \label{lem:contractionreduce}
\begin{enumerate}
\item
Any reduced word for a contraction element is a contraction word.
\item
Any non-trivial power of a contraction element is a contraction element.
\end{enumerate}
\end{lemma}

\proof
(1)
Let $w=\prod_{i=1}^p b_p^{e_i}$ be 
a contraction word of an anticonnected set $B$ in $V(\Gamma)$. 
Here, $b_i\in B$ and $e_i=\pm1$ for each $i$.
Suppose $w'$ is a reduced word, such that $w'=_{A(\Gamma)}w$.
There exists a 
 dual van Kampen diagram $\Delta$ for
$ww'^{-1}$. Note that any properly embedded arc of $\Delta$
 meets both of the intervals $w$ and $w'$,
since $w$ and $w'$ are reduced (Lemma~\ref{lem:qpair}). 
Now let $b,b'\in B$. 
$w$ is a contraction word, so one can find
$l\ge1$ and $1\le k_1< k_2 <\cdots<k_l\le p$ such that,
$b_{k_i}$ and $b_{k_{i+1}}$ are non-adjacent for each $i=1,\ldots,l-1$,
and $b=b_{k_1},b'=b_{k_l}$. Let $\gamma_i$ be the arc that intersects
with the segment $b_{k_i}$ of $w$. Since $\gamma_1,\gamma_2,\ldots,\gamma_l$
 are all disjoint,
the boundary points of those arcs on $w'$ will yield the desired subsequence
of the letters of $w'$.

(2) Let $u\pmo v u$ be a reduced word for $g$ as in 
Lemma~\ref{lem:conjugatepower}.
Note that a sequence, containing a contraction sequence as a monotonic
 subsequence, is again a contraction sequence.
So 
the reduced word
$u\pmo v^m u$ is a contraction word of $B$, for each $m\ne 0$.
\qed

\renewcommand{\labelenumi}{(\roman{enumi})}

\begin{definition} \label{def:canonical}
Let $\Gamma$ be a graph, and $P$ and $Q$ be disjoint
subsets of $V(\Gamma)$. Suppose $P_1$ is a set of words in
$\lll P\rrr\le A(\Gamma)$. 
A {\em  canonical expression for $g\in\lll P_1,Q\rrr$
with respect to $\{ P_1,Q\}$
} is a word
$\prod_{i=1}^k  c_i^{e_i}$, where
\begin{enumerate}
\item
$ c_i\in P_1\cup Q$
\item
$e_i=1$ or $-1$
\item
$\prod_{i=1}^k  c_i^{e_i}  =_{A(\Gamma)} g$
\end{enumerate}
such that $k$ is minimal. $k$ is called the {\em length}
of the canonical expression.
\end{definition}

\begin{remark}
In the above definition, a canonical 
expression exists for any element in $\lll P_1,Q\rrr$. 
In the case when $P_1\subseteq P$,
a word is a canonical expression with respect to 
$\{P_1,Q\}$,
if and only if it is reduced
in $A(\Gamma)$.
\end{remark}

Now we  compute
intersections of certain subgroups of $A(\Gamma)$.

\renewcommand{\labelenumi}{(\arabic{enumi})}
\begin{lemma} \label{lem:intersections}
Let $\Gamma$ be a graph, $P,Q$ be disjoint subsets of
$V(\Gamma)$ and
$P_1$ be a set of words in $\lll P\rrr\le A(\Gamma)$.
Let $R$ be any subset of $V(\Gamma)$.
\begin{enumerate}
\item
If $w$ is a canonical expression with respect to
$\{P_1,Q\}$, then
there does not exist a $q$-pair of $w$
for any $q\in Q$.
\item
$\lll P_1,Q\rrr\cap\lll R\rrr \subseteq\lll P_1,Q\cap R\rrr$.
Moreover, the equality holds if $P\subseteq R$.
\item
Let $\tilde w$ be a contraction word of $P$,
and $P_1=\{\tilde w\}$.
Assume $P\not\subseteq R$.
Then 
$\lll P_1,Q\rrr\cap\lll R\rrr =\lll Q\cap R\rrr$.
\end{enumerate}
\end{lemma} 

\proof

(1) 
Let $w$
be a canonical expression,
Suppose there exists
a $q$-pair of $w$ for some $q\in Q$. Then by Lemma~\ref{lem:qpair},
 one can 
write
$w=w_1 q^{\pm1}w_2q^{\mp1}w_3$ for some subwords
$w_1,w_2$ and $w_3$ such that 
$w_2\in\lll \link_\Gamma(q) \rrr$.
It follows that 
$w=_{A(\Gamma)}w''
 =w_1 w_2 w_3$.
Since $P\cap Q=\varnothing$, $w_1,w_2$ and $w_3$ are also
canonical expressions with respect to $\{P_1,Q\}$.
This contradicts to the minimality of $k$.

(2)
Let $w$
 be a canonical expression   and
$w'=_{A(\Gamma)}w$ be a reduced word.
Consider a dual van Kampen diagram $\Delta$  
for $w w'\pmo$.

Suppose that there exists a $q$-segment in $w$, 
for
some $q\in Q$.
Then by (1), the $q$-segment should be joined, by a $q$-arc, 
to another $q$-segment
of $w'$.
Since
$w'$ is a reduced word representing an element in $\lll R\rrr$,
each segment of $w'$ 
is labeled by
$R^{\pm1}$
(Lemma~\ref{lem:subgraph} (2)).
Therefore,
$q\in Q\cap R$.

If $P\subseteq R$, then 
$
\lll P_1,Q\cap R\rrr
\subseteq
\lll P_1,Q\rrr\cap\lll R\rrr
$ is obvious.

(3) 
$\lll Q\cap R\rrr\subseteq \lll P_1,Q\rrr\cap\lll R\rrr$ is obvious.

To prove the converse,
suppose $w\in(\lll P_1,Q\rrr\cap \lll R\rrr)\setminus\lll Q\cap R\rrr$.
$w$ is chosen so that $w$ is a canonical expression with respect to
$\{P_1,Q\}$, and the length (as a canonical expression) is minimal.

\begin{figure}[ht!]
\psfrag{w}{\footnotesize $w$}
\psfrag{wp}{\footnotesize $w'$}
\psfrag{w2}{\footnotesize $w_2$}
\psfrag{w1}{\footnotesize $w_1$}
\psfrag{wt}{\footnotesize $\tilde w$}
\psfrag{we}{\footnotesize $\tilde w^e$}
\psfrag{b}{\footnotesize $b$}
\psfrag{bm}{\footnotesize $b^{\mp1}$}
\psfrag{bp}{\footnotesize $b'$}
\psfrag{bpm}{\footnotesize $b'^{\mp1}$}
\psfrag{ds}{\footnotesize $\cdots$}
\psfrag{vd}{\footnotesize $\vdots$}
\psfrag{bt}{\footnotesize $\beta$}
\psfrag{b1}{\footnotesize $\beta_1$}
\psfrag{b2}{\footnotesize $\beta_2$}
\psfrag{bl}{\footnotesize $\beta_l$}
\centerline{\epsfig{file=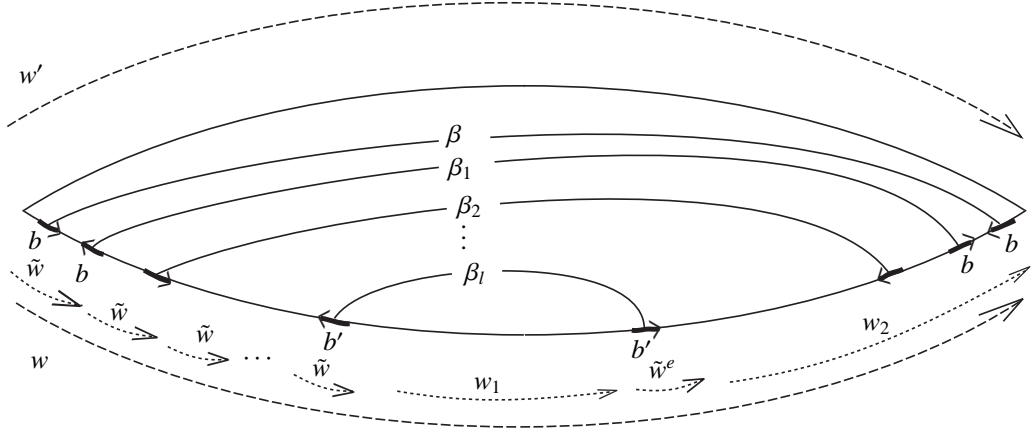, width=.90\hsize}}
\caption{$\Delta$ in the proof of Lemma~\ref{lem:intersections}.
\label{fig:intersections}}
\end{figure}

Let  $w=\prod_{i=1}^k c_i^{e_i}$ ($c_i\in\{P_1,Q\}$, $e_i=\pm1$),
$w'$ be a reduced word satisyfing $w'=_{A(\Gamma)}w$,
 and 
  $\Delta=(\HH,\lambda)$ be
a  
dual van Kampen diagram for $w w'\pmo$ (Figure~\ref{fig:intersections}).
Any shorter canonical expression than $w$, for an element in
$\lll P_1,Q\rrr\cap\lll R\rrr$ is in
$\lll Q\cap R\rrr$. This implies that $c_1=\tilde w=c_k$.
Note that
each segment of $w'$
is labeled by $R^{\pm1}$.
From the proof of (2),
$c_i\in P_1\cup (Q\cap R)=\{\tilde w\}\cup (Q\cap R)$ for each $i$. 

Now suppose $ c_i=\tilde w$ for some $i$. 
Fix  $b\in P\setminus R$. 
Choose the  $b$-arc $\beta$ that intersects with
the leftmost $b$-segment in $w$ on 
$\partial\Delta$. This $b$-segment is contained
in the leftmost $\tilde w$-interval in $w$.
By the minimality of the length of $w$, we can write
$w=\tilde w^{m} w_1 \tilde w^{e} w_2$ for some
subwords $w_1,w_2$ of $w$, $m\in\Z\setminus\{0\}$ and $e\in\{1,-1\}$. 
$w_1$ and $w_2$ are chosen so that
the letters of $w_1$ are in $(Q\cap R)^{\pm1}$
and $\beta$ intersects with a segment in the interval $\tilde w^{e} w_2$.
Without loss of generality,
we may assume $m>0$. 

Let $b'$ be any element in $P$. 
By Lemma~\ref{lem:contractionreduce}, any reduced word for $\tilde w^m$
is a contraction word of $P$. So, one can find a sequence of
arcs $\beta_1,\beta_2,\ldots,\beta_l\in\HH$ such that
\renewcommand{\labelenumi}{(\roman{enumi})}
\begin{enumerate}
\item
$\lambda(\beta_1)=b,\lambda(\beta_l)=b'$,
\item
$\lambda(\beta_i)$ and $\lambda(\beta_{i+1})$ are non-adjacent in $\Gamma$,
for each $i=1,2,\ldots,l-1$, and
\item
each $\beta_i$ intersects with a segment in the interval 
$\tilde w^{e}w_2$.
\end{enumerate}
Note that (iii) comes from the assumptions that $\beta_i$ does not
join two segments from $\tilde w^m$ (by reducing $\tilde w^m$ first),
and that the letters of $w_1$
are in $(Q\cap R)^{\pm1}$.

As in (2), each segment of $w_1$ is joined
to a segment in $w'$. In particular,
$[b',w_1]=
[\lambda(\beta_l),w_1]\eqag1$. Since this 
is true for any $b'\in P$, 
$w\eqag w_1 \tilde w^{m+e} w_2$.
One has
$\tilde w^{m+e} w_2\in 
(\lll P_1,Q\rrr\cap \lll R\rrr)\setminus
\lll Q\cap R\rrr$,
for
$w\not\in\lll Q\cap R\rrr$ and 
$w_1\in\lll Q\cap R\rrr$.
By the minimality of $w$, we have $w_1=1$.
This argument continues, and finally one can write
 $w=\tilde w^{m'}$ for some $m'\ne0$. In particular,
 any reduced word for $w$ is a contraction word of $P$
  (Lemma~\ref{lem:contractionreduce}).
This is impossible since $w\in\lll R\rrr$ and $P\not\subseteq R$.
\qed

\begin{lemma} \label{lem:contraction_element}
Let $\Gamma$ be a graph, $B$ be an anticonnected set of vertices of $\Gamma$
and $g$ be a contraction element of $B$.
Then 
there exists an injective map
$ \phi:A (\cob(\Gamma,B))\rightarrow A(\Gamma)$
satisfying 
\[\phi(x) =
\left\{
\begin{array}{ll}
g\quad & \textrm{if }x = v_{B}\\
x\quad& \textrm{if }x\in V(\Gamma)\setminus B
\end{array}
\right.
\]
\end{lemma}

\proof

As in the proof of 
Lemma~\ref{lem:contract2}, let $\hat\Gamma=\cob(\Gamma,B),\hat v=v_B$
and $A=\{q:q\in V(\Gamma)\setminus B\}$.
For $q\in A$,
let $\hat q$ denote the corresponding vertex in $\hat\Gamma$, and 
$\hat A=\{\hat q:q\in A\}$. There exists a map 
$\phi:A(\hat\Gamma)\rightarrow A(\Gamma)$
satisfying
\[\phi(x) =
\left\{
\begin{array}{ll}
g &\textrm{if } x = \hat v\\
q &\textrm{if } x=\hat q\in\hat A
\end{array}
\right.
\]
To prove that $\phi$ is injective, we use an induction on $|A|$.

If $A=\varnothing$, then $V(\Gamma)=B$ and
$\hat\Gamma$ is the graph with one vertex $\hat v$. 
So, $\phi$ maps $\lll \hat v\rrr  =A(\hat\Gamma)\cong\Z$ isomorphically
onto $\Z\cong\lll g\rrr\le A(\Gamma)$. 

Assume the injectivity of $\phi$ 
for the case when $|A|=k$, and now let $|A|=k+1$.

Choose any $t\in A$. Let $A_0=A\setminus\{t\}$ and 
$\hat A_0=\{\hat q:q\in A_0\}$.
Let 
$\Gamma_0$ be the induced subgraph on $A_0\cup B$ in $\Gamma$,
and
$\hat\Gamma_0$ be the induced subgraph on $\hat A_0\cup\{\hat v\}$ 
in $\hat\Gamma$.
We consider $A(\Gamma_0)$ and $A(\hat\Gamma_0)$ as subgroups of
 $A(\Gamma)$ and $A(\hat\Gamma)$, respectively, so that
$A(\Gamma_0)=\lll A_0,B\rrr$
and
 $A(\hat\Gamma_0)=\lll\hat A_0,\hat v\rrr$.
Let $K=\lll A_0,g\rrr=\phi(A(\hat\Gamma_0))$ and 
$J=\lll A,g\rrr=\phi(A(\hat\Gamma))$.
By the inductive hypothesis, $\phi$ 
maps $A(\hat\Gamma_0)$ isomorphically
onto $K$ (Figure~\ref{fig:contract_element}).

\begin{figure}[ht!]
\[
\xymatrix{
&&& A(\Gamma)= A(\Gamma_0)\ast_C\ar@{-}[d]\ar@{-}[dl]\\
A(\hat\Gamma)=A(\hat\Gamma_0)\ast_{D}
\ar@{-}[d]\ar[rr]^(.6)*{\phi} && J=\lll A, 
g\rrr\ar@{-}[d]
& A(\Gamma_0)=\lll A_0,B\rrr\ar@{-}[dl]\ar@{-}[d]\\
\quad A(\hat\Gamma_0)=\lll\hat A_0,\hat v\rrr\quad 
\ar[rr]^(.6)*{\phi\downharpoonright}_(.6)*{\cong}\ar@{-}[d]
&&
K=\lll A_0,g\rrr
& C=\lll\link_\Gamma(t)\rrr\ar@{-}[d]\\
D = \lll\link_{\hat\Gamma}(\hat t)\rrr\ar@{.>}[rrr]&&&
J\cap C
}\]
\caption{Proof of Lemma~\ref{lem:contraction_element}.
Note that $V(\Gamma)=A\sqcup B=A_0\cup\{t\}\cup B$ and
$V(\hat\Gamma)=\hat A\sqcup \{\hat v\}=\hat A_0\cup\{\hat t\}\cup\{\hat v\}$.
\label{fig:contract_element}}
\end{figure}

From Lemma~\ref{lem:artinhnn}, we can identify
$A(\Gamma) = A(\Gamma_0)\ast_C$,  where
$C=\lll \link_\Gamma(t)\rrr$
and
$t$ is the stable letter.
 Since $J=\lll A_0,g,t\rrr=\lll K,t\rrr$,
 Lemma~\ref{lem:jkt} implies that we can also
 identify
$J=(J\cap A(\Gamma_0))\ast_{J\cap C}$, where $t$ is
the stable letter again.
Also, we identify
$A(\hat\Gamma) = A(\hat\Gamma_0)\ast_{D}$,
where 
$D  = \lll\link_{\hat\Gamma}(\hat t)\rrr$
and
$\hat t$ is the stable letter.

By Lemma~\ref{lem:intersections} (2),
$J\cap A(\Gamma_0) =\lll g,A\rrr\cap\lll A_0,B\rrr
=\lll g,A\cap(A_0\cup B)\rrr
=\lll g,A_0\rrr =\phi(A(\hat\Gamma_0))$.

Applying Lemma~\ref{lem:intersections} (3) for the case when
$R=\link_\Gamma(t)$,
\begin{eqnarray*} 
J\cap C &=& 
\lll g,A\rrr\cap\lll \link_\Gamma(t)\rrr\\
&=&
\left\{
\begin{array}{cl}
\lll\link_\Gamma(t)\cap A,
g \rrr &\quad \textrm{if }B\subseteq \link_\Gamma(t) \\
\lll\link_\Gamma(t)\cap A\rrr &\quad\textrm{otherwise }
\end{array}
\right.\end{eqnarray*}
From the definition of a co-contraction, we note that
\[
D = \link_{\hat\Gamma}(\hat t)=
\left\{
\begin{array}{ll}
\{\hat q:q\in\link_\Gamma(t)\cap A\}\cup\{\hat v\}
 &\textrm{if } B\subseteq \link_\Gamma(t) 
 \\
 \{\hat q:q\in\link_\Gamma(t)\cap A\}
 &\textrm{otherwise }
\end{array}
\right.
\]
Hence, $J\cap C = \phi(D)$. This implies that
                     $\phi:A(\hat\Gamma)\rightarrow
                     J$ is an isomorphism, as follows.
\[
\xymatrix{
D
\ar[d]^*{\cong} & 
\le & A(\hat\Gamma_0)\ar[d]^*{\cong} &\le
& A(\hat\Gamma_0)\ast_{D}
& = & A(\hat \Gamma)\ar[d]_*{\phi}\\
J\cap C&\le& K=J\cap A(\Gamma_0) & \le & (J\cap A(\Gamma_0))\ast_{J\cap C}
&=& J
}
\]
\qed

Now the following theorem is immediate by an induction on $m$.

\begin{theorem} \label{thm:anticntgen}
Let $\Gamma$ be a graph and $B_1$,$B_2$,\ldots,$B_m$ be
disjoint subsets of $V(\Gamma)$ such that each $B_i$ is anticonnected.
For each $i$, let $v_{B_i}$ denote the vertex correponding to 
$B_i$ in $\cob(\Gamma,(B_1,B_2,\ldots,B_m))$, and
$g_i$ be a contraction element of $B_i$.
Then 
there exists an injective map
$ \phi:A (\cob(\Gamma,(B_1,B_2,\ldots,B_m) ))\rightarrow A(\Gamma)$
satisfying 
\[\phi(x) =
\left\{
\begin{array}{ll}
g_i\quad &\textrm{if } x = v_{B_i}\ ,\ \textrm{for some }$i$\\
x\quad&\textrm{if } x\in V(\Gamma)\setminus\cup_{i=1}^m B_i\\
\end{array}
\right.
\]\qedhere
\end{theorem}

We conclude this article by noting that 
there is 
another partial answer to 
the question of which right-angled Artin groups
contain hyperbolic surface groups. 
Namely, if $\Gamma$ does not contain an induced cycle of length $\ge5$,
and either $\Gamma$ does not contain an induced $C_4$ (hence chordal),
or $\Gamma$ is triangle-free (hence bipartite), then
$A(\Gamma)$ does not contain a hyperbolic surface group~\cite{kim2006}.

\bibliographystyle{amsplain}
\bibliography{cb2006} 

\end{document}